\documentclass[reqno,english,12pt]{amsart}

\pdfoutput=1                             

\usepackage{float}
\usepackage[latin1]{inputenc}            
\usepackage[T1]{fontenc}                 
\usepackage[english]{babel}      
\usepackage{color}
\usepackage{fullpage}
\usepackage{verbatim}

\usepackage{mathtools}
\usepackage{euscript}
\usepackage{url}

\newtheorem{theorem}{Theorem}

\newtheorem{definition}[theorem]{Definition}
\newtheorem{lemma}[theorem]{Lemma}
\newtheorem{proposition}[theorem]{Proposition}
\newtheorem{corollary}[theorem]{Corollary}

\newcommand{\R}{\mathbb{R}}
\newcommand{\N}{\mathbb{N}}

\newcommand{\Z}{\mathbb{Z}}
\newcommand{\Q}{\mathbb{Q}}

\def\co{\colon\thinspace}

\usepackage{amssymb}
\usepackage{yfonts}
\def\QED{\hfill{$\Box$}}
\def\proof{\textit{Proof.} }
\def\Ker{\text{Ker} }
\def\H{\mathbb{H}^2 }
\def\lan{\langle}
\def\ran{\rangle}

\def\O{\mathcal{O}}
\def\Hom{\text{Hom}}

\begin{document}

\title{Zariski dense surface subgroups in $SL(n,\Q)$ with odd $n$}
\author{Carmen Galaz-Garc\'ia}
\begin{abstract}
For odd $n$ we construct a path $\rho_t\co \pi_1(S) \to SL(n,\R)$ of discrete, faithful and Zariski dense representations of a surface group such that $\rho_t(\pi_1(S)) \subset SL(n,\Q)$ for every $t\in \Q$. 

\end{abstract}

\maketitle


Constructing Zariski dense surface subgroups in $SL(n,\R)$ has  attracted attention as a step to finding \textit{thin groups}, these are infinite index subgroups of a lattice in $SL(n,\R)$ which are Zariski dense. 
Finding thin subgroups inside lattices in a variety of Lie groups has been a topic of significant interest in recent years, in part from the connections thin groups have to expanders and the affine sieve of Bourgain, Gamburd, and Sarnak \cite{bourgain2010affine}\cite{sarnak2014notes}. 

Though thin subgroups are in a sense generic \cite{Fuchs_14_ubiquity}\cite{fuchs2017generic}, finding particular specimens of thin surface subgroups in a given lattice remains a difficult task.
In this direction Long, Reid and Thistlethwaite \cite{Long_11_SL3Z} produced in 2011 the first infinite family of nonconjugate thin surface groups in $SL(3,\Z)$. 
Their approach relies on parametrizing a family of representations $\rho_t$ of the triangle group $\Delta(3,3,4)$ in the \textit{Hitchin component}, so that for every $t\in \Z$ the subgroup $\rho_t(\Delta(3,3,4))$ is in $SL(3,\Q)$ and has integral traces. 
By results of Bass \cite{bass1980groups} these two properties together with $\rho_t(\Delta(3,3,4))$ being non-solvable and finitely generated guarantee that it is conjugate to a subgroup of $SL(3,\Z)$. 
In 2018 Long and Thistlethwaite \cite{Long_18_SL4Z} used a similar approach  to obtain an infinite family of non-conjugate Zariski dense surface subgroups in $SL(4,\Z)$ and $SL(5,\Z)$. 

Ballas and Long \cite{Ballas_18_bending} in turn used the idea of "bending" a representation of the fundamental group of a hyperbolic $n$-manifold $\pi_1(N)$ along an embedded totally geodesic and separating  hypersurface to obtain thin groups in $SL(n+1,\R)$ which are isomorphic to $\pi_1(N)$. 
The goal of this article is to combine the aforementioned approaches to construct a family of Zariski dense rational surface group representations by bending orbifold representations. 
Our main result is the following:
\begin{theorem}\label{thm1}
For every surface $S$ finitely covering the orbifold $\O_{3,3,3,3}$ and every odd $n>1$ there exists a path of discrete, faithful and irreducible representations $\rho_t \co \pi_1(S) \to SL(n,\R)$, so that
\begin{enumerate}
	\item  $\rho_0(\pi_1(S)) < SL(n,\Z)$, 
	\item $\rho_t$ is Zariski dense for every $t>0$, and
	\item $\rho_t(\pi_1(S)) < SL(n,\Q)$ for every $t\in \Q$. 
\end{enumerate}
\end{theorem} 

Every representation $\rho_t$ in theorem \ref{thm1} is a surface Hitchin representation.  
Several of its properties are derived from the seminal work of Labourie \cite{Labourie_06_AnosovReps} on Anosov representations, the classification of Zariski closures of surface Hitchin representations by Guichard \cite{Guichard_ZariskiClosure} and  the recent introduction of orbifold Hitchin representations by Alessandrini, Lee and Schaffhauser \cite{ALS_19_orbifoldHitchin}. 
We provide an overview of these results in sections \ref{sectionHitchin} and \ref{ZdenseSection}. 
At the end of section \ref{ZdenseSection} we also prove the following criterion for Zariski density, which will be subsequently used to discard Zariski closures.

 \begin{proposition}\label{prop1.2}
 Let $\rho \co \pi_1(\O) \to PSL(n,\R)$ be an orbifold Hitchin representation such that 
\begin{itemize}
\item if $n=2k$ is even then $\rho(\pi_1(\O))$ is not conjugate to a subgroup of $PSp(2k,\R)$ or,
\item if $n=2k+1$ is odd then  $\rho(\pi_1(\O))$ is not conjugate to a subgroup of $PSO(k,k+1)$.
\end{itemize}
Then $\rho(H)$ is Zariski dense in $PSL(n,\R)$ for every finite index subgroup $H$ of $\pi_1(\O)$. 
\end{proposition}
In section \ref{sectionBending} we give a general construction to obtain a path of representations as in theorem \ref{thm1}.
This is based on bending  the fundamental group $\pi_1(\O)$ of a hyperbolic 2-dimensional orbifold along a simple closed curve in $\O$ with infinite order as an element of $\pi_1(\O)$. 
Theorem \ref{thm1} then follows from applying the results in section 2 to a suitable representation of the fundamental group of the orbifold $\O_{3,3,3,3}$ whose underlying topological space is $S^2$ and has four cone points of order 3. 
This final step is covered in section \ref{sectionO3333}.

\textbf{Remark.} During the finalization of this project, Long and Thistlethwaite used bending to construct thin surface groups in $SL(n,\Z)$ for every odd $n$ \cite{Long_This_20_SLodd}, the even case remains open. 


\section{Hitchin representations}\label{sectionHitchin}

In this section we give a short overview of surface and orbifold Hitchin representations. 

\subsection{Spaces of representations} \label{section_irrep} 

\sloppy Let $G$ be a Lie group and let $\Gamma$ be a group with a finite presentation $\lan \alpha_1, \ldots, \alpha_k \ | \ r_1, \ldots , r_m \ran$. 
Then every relator $r_i$ defines a map $R_i \co G^k \to G$. 
If we let $\Hom(\Gamma, G) = \cap_{i=1}^m R_i^{-1}(Id)$, then the map $\phi \mapsto (\phi(\alpha_1), \ldots, \phi(\alpha_k))$ is a bijection between the set of all group homomorphisms from $\Gamma$ to $G$ and $\Hom(\Gamma, G)$. 
We will regard $\Hom(\Gamma, G)$ as having the subspace topology from $G^k$. 

Let $\Hom^+(\Gamma, G)$ be the subset of representations in $\Hom(\Gamma,G)$ which decompose as a direct sum of irreducible representations and let
	$\text{Rep}^+(\Gamma, G) = \Hom^+(\Gamma,G)/G$
be the quotient space by the conjugation action. 
With the quotient topology $\text{Rep}^+(\Gamma, G)$ has the structure of an algebraic variety (\cite{BGGW_07_RepsNotes} sec. 5.2).

In the following we will  frequently use the representation 
\begin{eqnarray}\label{SLirrep}
	\tilde\omega_n \co SL(2,\R) \to SL(n,\R)
\end{eqnarray}
 given by the action of $SL(2,\R)$ on the vector space $\mathcal{P}$ of homogeneous polynomials in 2 variables of degree $n-1$. 
If $n = 2k$ is even, the image of $\tilde\omega_n$ is contained in the symplectic group $Sp(2k, \R)$, and if $n = 2k + 1$ is odd, it is contained in a group isomorphic to $SO(k, k + 1)$.

It is well known that the representation $\tilde\omega_n$ is absolutely irreducible and is, up to conjugation, the unique irreducible representation from $SL(2,\R)$ into $SL(n,\R)$.
This representation induces a \textit{ projective representation} $\omega_n \co PSL(2,\R) \to PSL(n,\R)$ which is also irreducible and unique up to conjugation.

\subsection{Hitchin representations of surface groups} 

Let $S$ be a closed surface of genus $g>1$.  
In 1988 Goldman proved that $\text{Rep}^+(\pi_1(S),PSL(2,\R))$ has $4g-3$ connected components, 
two of which are diffeomorphic to $\R^{6g-6}$ and called these \textit{Teichm\"uller spaces} (\cite{Goldman_88_TopComponents} thm. A, see also note at end of  thm. 10.2 in \cite{Hitchin92}).  
The two Teichm\"uller spaces $\mathcal{T}^\pm(S)$ are precisely the sets of conjugacy classes by $PSL(2,\R)$ of \textit{Fuchsian representations}, which are discrete and faithful representations $\rho\co \pi_1(S) \to PSL(2,\R)\equiv\text{Isom}^+(\H)$. 

\begin{definition}\label{fuchsian_def}\em
For $n>2$ a representation $r\co \pi_1(S) \to PSL(n,\R)$ is called \textit{Fuchsian} if it can be decomposed as $r = \omega_n\circ r_0$ where $r_0\co \pi_1(S)\to PSL(2,\R)$ is discrete and faithful, and $\omega_n\co PSL(2,\R)\to PSL(n,\R)$ is the unique irreducible representation introduced in section \ref{section_irrep}.
\end{definition}
\begin{definition}\em
The \textit{Fuchsian locus} is the set of all $PSL(n,\R)$ conjugacy classes of Fuchsian representations, namely the set $\omega_n(\mathcal{T}^\pm(S))$.
\end{definition}
The space $\text{Rep}^+(\pi_1(S), PSL(n,\R))$ has three topological connected components if $n$ is odd and 6 if $n$ is even  (\cite{Hitchin92}, thm. 10.2). 
The Fuchsian locus is contained in one component in the odd case and in two components in the even case. 
Each of these distinguished components, called \textit{Hitchin components}, is diffeomorphic to $\R^{(1-n^2)(1-g)}$. 
When $n>2$ is even, both Hitchin components are related by an inner automorphism of $PSL(n,\R)$. 
In the odd case, where there is only one component, we will denote the Hitchin component by $\text{Hit}(\pi_1(S),PSL(n,\R))$.

\begin{definition}\em
Let $S$ be a closed surface of genus greater than one.
A representation $r \co \pi_1(S) \to PSL(n,\R)$ is a \textit{surface Hitchin representation} if its $PSL(n,\R)$-conjugacy class belongs to a Hitchin component of $\text{Rep}^+(\pi_1(S),PSL(n,\R))$. 
\end{definition}

In \cite{Labourie_06_AnosovReps}, Labourie introduces \textit{Anosov representations} and proves that surface Hitchin representations are $B$-Anosov where $B$ is any Borel subgroup of $PSL(n,\R)$. 
This gives surface Hitchin representations essential algebraic properties, out of which we will use theorem \ref{Hitchin_reps} below.

\begin{definition}[\cite{bridgeman2020simple} sec. 2.2]\label{loxo_def}\em
A matrix $A\in SL(n,\R)$ is \textit{purely loxodromic} if it is diagonalizable over $\R$ with eigenvalues of distinct modulus. If $A \in PSL(n,\R)$ then we say $A$ is \textit{purely loxodromic} if any lift of $A$ to an element of $SL(n,\R)$ is purely loxodromic. 
\end{definition}

\begin{theorem}[\cite{Labourie_06_AnosovReps} thm. 1.5, lemma 10.1]\label{Hitchin_reps}
 A surface Hitchin representation $r\co \pi_1(S) \to PSL(n,\R)$ is discrete, faithful and strongly irreducible. 
Moreover, the image of every non-trivial element of $\pi_1(S)$ under $r$ is purely loxodromic. 
 \end{theorem}
 
\subsection{Hitchin representations of orbifold groups} 

Let $\O$ be a 2-dimensional closed orbifold of negative orbifold Euler characteristic $\chi(\O)$ and let $\pi_1(\O)$ be its orbifold fundamental group.  
In \cite{Thurston_78_GeometryNotes} Thurston proves there is a connected component of the representation space $\text{Rep}(\pi_1(\O), PGL(2,\R))$
that parametrizes hyperbolic structures on $\O$.
This component is called the \textit{Teichm\"uller space} of the orbifold $\O$, we will denote it by $\mathcal{T}(\O)$. 
As with surfaces, the orbifold Teichm\"uller space consists of conjugacy classes of discrete and faithful representations of $\pi_1(\O)$ into $PGL(2,\R) \equiv \text{Isom}(\H)$, which we will call \textit{Fuchsian representations} too. 
More recently, Alessandrini, Lee, and Schaffhauser used the irreducible representation $\omega_n$   to define the \textit{Hitchin component} $\text{Hit}(\pi_1(\O),PGL(n,\R))$ of $\text{Rep}(\pi_1(\O),PGL(n,\R))$ as the unique connected component in this representation space which contains the connected Fuchsian locus $\omega_n(\mathcal{T}(\O))$ (\cite{ALS_19_orbifoldHitchin} def. 2.3) and prove $\text{Hit}(\pi_1(\O),PGL(n,\R))$ is homeomorphic to an open ball (\cite{ALS_19_orbifoldHitchin} thm. 1.2).

\begin{definition}[\cite{ALS_19_orbifoldHitchin} def. 2.4]\em
Let $\O$ be a 2-dimensional  connected closed orbifold with negative orbifold Euler characteristic. 
A representation $r\co \pi_1(\O) \to PGL(n,\R)$ is an \textit{orbifold Hitchin representation} if its $PGL(n,\R)$-conjugacy class belongs to the Hitchin component $\text{Hit}(\pi_1(\O), PGL(n,\R))$ of $\text{Rep}(\pi_1(\O),PGL(n,\R))$.
\end{definition}

The definition of Anosov representations has been generalized by Guichard and Wienhard (\cite{GuichardWienhard_12_Anosov} def. 2.10) to include representations of word hyperbolic groups into semisimple Lie groups. 
With this more general definition, and just as their surface counterparts, orbifold Hitchin representations are also $B$-Anosov where $B$ is a Borel subgroup of $PGL(n,\R)$ (\cite{ALS_19_orbifoldHitchin} prop. 2.16) and thus share some strong algebraic properties. 

\begin{theorem}[\cite{ALS_19_orbifoldHitchin} thm. 1.1]\label{orbifoldHitchin}
An orbifold Hitchin representation $r\co \pi_1(\O) \to PGL(n,\R)$ is discrete, faithful and strongly irreducible. 
Moreover, the image of every infinite order element of $\pi_1(\O)$ under $r$ is purely loxodromic. 
\end{theorem}

\section{Zariski dense Hitchin representations}\label{ZdenseSection}

In this section we focus on Zariski density of Hitchin representations and prove corollary \ref{corollary_Zdense} which gives a criterion to determine when the image of a finite index subgroup of an orbifold group under a Hitchin representation is Zariski dense. 

\subsection{Zariski closures of Hitchin representations}\label{subsec_irrep}
Let $G$ be an algebraic matrix Lie group, then $G$ has both its standard topology as a subset of some $\R^N$ and the Zariski topology. 
If $X$ is a subset of $G$ then its \textit{Zariski closure} is the closure of $X$ in $G$ with respect to the Zariski topology. 
We say a subgroup $H<G$ is \textit{Zariski dense} in $G$ if its Zariski closure equals $G$.  
A representation $r\co \Gamma \to G$ is \textit{Zariski dense} if $r(\Gamma)$ is Zariski dense in $G$. 

The image of the irreducible representation $\omega_n\co PSL(2,\R) \to PSL(n,\R)$ is contained, if $n$ is even, in a conjugate of $PSp(n, \R)$, which is the projectivization of the symplectic group $Sp(n, \R)$.
If $n = 2k + 1$ is odd, the image of $\omega_n$ is contained in a conjugate of the orthogonal group  $SO(k, k + 1) = PSO(k,k+1)$.
This implies that the images of Fuchsian representations are contained in (a conjugate of) $PSp(n,\R)$ or in $SO(k,k+1)$ and, in particular, they are not Zariski dense. 
More generally, for surface Hitchin representations Guichard \cite{Guichard_ZariskiClosure} has  announced a classification of Zariski closures of their lifts.
An alternative proof of this result has been given recently by Sambarino (\cite{Sambarino_closures} cor. 1.5). 
The version of this result we cite here comes from theorem 11.7 in \cite{BCLS_15_PressureMetric}.

\begin{theorem}[\cite{Guichard_ZariskiClosure}, \cite{Sambarino_closures}] \label{GuichardResult}
If $r \co \pi_1(S) \to SL(n,\R)$ is the lift of a surface Hitchin representation and $H$ is the Zariski closure of $r(\pi_1(S))$, then 
\begin{itemize}
\item If $n=2k$ is even, $H$ is conjugate to either $\omega_n(SL(2,\R))$, $Sp(2k,\R)$ or $SL(2k,\R)$. 
\item If $n=2k+1$ is odd and $n\neq7$, then $H$ is conjugate to either $\omega_n(SL(2,\R))$, $SO(k,k+1)$ or $SL(2k+1,\R)$. 
\item If $n=7$, then $H$ is conjugate to either $\omega_7(SL(2,\R))$, $G_2$, $SO(3,4)$ or $SL(7,\R)$.
\end{itemize}
\end{theorem}

\subsection{A criterion for Zariski density} 

Here we prove proposition \ref{prop1.2} which gives us a criterion to find Zariski dense Hitchin representations. 

\begin{lemma}\label{lift_even}
Let $\rho \co \pi_1(\O) \to PSL(n,\R)$ with $n$ even be an orbifold Hitchin representation. 
Then for every $[\alpha] \in \pi_1(\O)$ of infinite order there is a lift $A \in SL(n,\R)$ of $\rho([\alpha])$ which has $n$ positive distinct eigenvalues.
\end{lemma}
\proof 
First consider a Fuchsian representation $\sigma \co \pi_1(\O) \to PSL(2,\R)$ and  $[\alpha]$ an infinite order element of $\pi_1(\O)$. 
Since $\O$ is a hyperbolic orbifold, $\sigma([\alpha])$ is conjugate to a hyperbolic element
$\begin{bmatrix} \lambda & 0 \\ 0 & \frac{1}{\lambda} \end{bmatrix} \in PSL(2,\R)$. 
We can lift this element to a matrix 
$\begin{pmatrix} \lambda & 0 \\ 0 & \frac{1}{\lambda} \end{pmatrix} \in SL(2,\R)$ with $\lambda >0$. 
Let $\tilde\omega_n \co SL(2,\R) \to SL(n,\R)$ be the unique irreducible representation in (\ref{SLirrep}), 
then $\tilde\omega_n\begin{pmatrix} \lambda & 0 \\ 0 & \frac{1}{\lambda} \end{pmatrix} \in SL(n,\R)$ has $n$ distinct positive eigenvalues $\lambda^{n-1}, \lambda^{n-3}, \ldots, \lambda^{-(n-3)}, \lambda^{-(n-1)}$ and is a lift of $\omega_n\circ\sigma([\alpha])\in PSL(n,\R)$. 

Now consider a Hitchin representation  $\rho\co \pi_1(\O) \to PSL(n,\R)$. 
Let $\rho_t$ be a path of Hitchin representations such that $\rho_0$ is Fuchsian and $\rho_1=\rho$. 
This induces a path $\rho_t([\alpha]) \subset PSL(n,\R)$. 
By the previous argument we may lift $\rho_t([\alpha])$ to a path $\tilde{A}_t \in SL(n,\R)$ such that $\tilde{A}_0$ has $n$ distinct positive eigenvalues. 
Since each eigenvalue of $\tilde{A}_t$ varies continuously and $\det\tilde{A}_t \neq 0$, all eigenvalues of $\tilde{A}_t$ are positive. 
Moreover, by theorem \ref{orbifoldHitchin} the absolute values of the eigenvalues of $\rho_t([\alpha])$ are distinct. 
This in turn implies all the eigenvalues of $\tilde{A}_t$ are distinct. 
Therefore $\tilde{A}_1 \in SL(n,\R)$ is a lift of $\rho([\alpha])$ with $n$ positive distinct eigenvalues. 

\QED

To prove our criterion for Zariski density (propositions \ref{Zdense_even} and \ref{Zdense})  we will make use of the following theorem by Culver. 

\begin{theorem}[\cite{Culver} thm. 2]\label{Culver_thm}
 Let $C$ be a real square matrix. Then the equation $C=\exp(X)$ has a unique real solution $X$ if and only if all the eigenvalues of $C$ are positive real and no elementary divisor (Jordan block) of C belonging to any eigenvalue appears more than once.
\end{theorem}

\begin{proposition}\label{Zdense_even}
Let $\rho\co \pi_1(\O) \to PSL(n,\R)$ with $n$ even be an orbifold Hitchin representation  so that $\rho(\pi_1(\O))$ is not conjugate to a subgroup of $PSp(n,\R)$. 
If $S$ is a surface finitely covering $\O$ then $\rho(\pi_1(S))$ is Zariski dense.
\end{proposition}
\proof 
Let $S$ be a surface finitely covering $\O$ and  
suppose that $\rho(\pi_1(S))$ is conjugate to a subgroup of $PSp(n,\R)$. 
Then there exists an alternating form $\Omega\in SL(n,\R)$ such that 
	$Sp(\Omega) = \{ g \in SL(n,\R) \ | \ g^T\Omega g = \Omega \}$
and $\rho(\pi_1(S)) \subset PSp(\Omega) = Sp(\Omega)/\pm I$.

Let $[\alpha] \in \pi_1(\O)$ be an infinite order element. 
By lemma \ref{lift_even} we can lift $\rho([\alpha]) \in PSL(n,\R)$ to a matrix $A\in SL(n,\R)$ with $n$ positive distinct eigenvalues. 
Since $\pi_1(S)$ has finite index in $\pi_1(\O)$ there exists a $k\in \N$ such that $\rho([\alpha])^k \in \rho(\pi_1(S))$. 
Then $A^k$ is a lift of $\rho([\alpha])^k$ and $A^k \in Sp(\Omega)$. 
Given that $A$ has $n$ positive distinct eigenvalues, by theorem \ref{Culver_thm} there is a unique $X\in M_{n\times n}(\R)$ such that $\exp(X) = A$. 
Then using that  $\exp(kX) = A^k$ preserves $\Omega$ we get that
\begin{eqnarray*}
	\exp(kX)^T\Omega \exp(kX) = \Omega &\Rightarrow& \Omega^{-1}\exp(kX)^T\Omega = \exp(kX)^{-1}\\
					  &\Rightarrow& \exp(\Omega^{-1}(kX)^T\Omega) = \Omega^{-1}\exp(kX)^T\Omega =\exp(-kX).
\end{eqnarray*}
Applying  theorem \ref{Culver_thm} now to  $\Omega^{-1}\exp(kX)^T\Omega$ we obtain that \begin{eqnarray*}
					  \Omega^{-1}(kX)^T\Omega = -kX
					  &\Rightarrow& -\Omega(kX)^T\Omega = -kX\\
					  &\Rightarrow& \Omega(kX)^T\Omega = kX.
\end{eqnarray*}
This implies that $kX \in \mathfrak{sp}(\Omega)$ and thus $A = \exp(X)\in Sp(\Omega)$. 
Given that $A$ is a lift of $\rho([\alpha])$, we have that $\rho([\alpha]) \in PSp(\Omega)$. 
Since $\pi_1(\O)$ is generated by its infinite order elements we get that $\rho(\pi_1(\O)) \subset PSp(\Omega)$, a contradiction.
So it cannot be that $\rho(\pi_1(S))$ is conjugate to a subgroup of $PSp(n,\R)$. 
In particular, if $r$ is a lift of the Hitchin surface representation $\rho|_{\pi_1(S)}$ then the Zariski closure of $r(\pi_1(S))$ cannot be conjugate to a subgroup of $Sp(n,\R)$. 
By theorem \ref{GuichardResult} it must be that the Zariski closure of $r(\pi_1(S))$ is $SL(n,\R)$. 
Therefore the Zariski closure of $\rho(\pi_1(S))$ is $PSL(n,\R)$. 

\QED

In the case when $n=2k+1$ is odd, by theorem \ref{GuichardResult} the Zariski closure of $\rho(\pi_1(S))$ where $\rho$ is a surface Hitchin representation is either conjugate to a subgroup of $SO(k,k+1)$ or equals $SL(n,\R)$. 
By assuming there exists a symmetric bilinear form $J$ such that $\rho(\pi_1(S)) \subset SO(J)$ we have an analogous proof to that of \ref{Zdense_even} to get a criterion for Zariski density of surface Hitchin representations in the odd case. 

\begin{proposition}\label{Zdense}
Let $\rho\co \pi_1(\O) \to SL(n,\R)$ with $n$ odd be an orbifold Hitchin representation such that there is no real quadratic form $J$ for which $\rho(\pi_1(\O))\subset SO(J)$. 
If $S$ is a surface finitely covering $\O$ then $\rho(\pi_1(S))$ is Zariski dense.
\end{proposition}

Given that any finite index subgroup of $\pi_1(\O)$ contains a surface subgroup which has finite index in $\pi_1(\O)$ we obtain the following result. 

\begin{proposition}\label{corollary_Zdense}
Let $\rho \co \pi_1(\O) \to PSL(n,\R)$ be an orbifold Hitchin representation such that 
\begin{itemize}
\item if $n=2k$ is even then $\rho(\pi_1(\O))$ is not conjugate to a subgroup of $PSp(2k,\R)$ or,
\item if $n=2k+1$ is odd then  $\rho(\pi_1(\O))$ is not conjugate to a subgroup of $PSO(k,k+1)$.
\end{itemize}
Then  for every finite index subgroup $H$ of $\pi_1(\O)$ the image $\rho(H)$ is Zariski dense in $PSL(n,\R)$. 
\end{proposition}

\section{Bending representations of orbifold groups }\label{sectionBending}

Theorem \ref{Zdense_surfacegrps} in this section gives a general construction of a path $\rho_t$ of Zariski dense Hitchin surface representations into $SL(n,\R)$ for odd $n$. 
By requiring that the initial representation $\rho_0$ has image inside $SL(n,\Q)$ we obtain corollary \ref{rational_reps}, in which every representation $\rho_t$ with $t\in \Q$ also has image in $SL(n,\Q)$. 

\subsection{Bending representations}
Let $\O$ be a 2-dimensional orientable connected closed orbifold of negative orbifold Euler characteristic and  $\O_L$, $\O_R$ be open connected suborbifolds with connected intersection $\O_L\cap \O_R$. 
Given a representation $\rho \co \pi_1(\O)  \to G$ there is a standard way of "bending" $\rho$ by an element $\delta$ of the centralizer in $G$ of $\rho(\pi_1(\O_L\cap \O_R))$ to obtain a representation $\rho_\delta \co \pi_1(\O) \simeq \pi_1(\O_L)\ast_{\pi_1(\O_L\cap\O_R)}\pi_1(\O_R)\to G$ so that 
	$\rho_\delta(\pi_1(\O)) = \lan \rho(\pi_1(\O_L)), \delta \rho(\pi_1(\O_R))\delta^{-1}\ran$ (see for example \cite{Goldman_87_geometricstructures} sec. 5).

From now onwards we will consider the case where there is a simple closed curve $\gamma \subset \O$, not parallel to a cone point, that divides $\O$ into two orbifolds $\O_L$ and $\O_R$ which share $\gamma$ as their common boundary, so that
	$\pi_1(\O) \simeq \pi_1(\O_L)\ast_{\lan [\gamma] \ran}\pi_1(\O_R).$

\begin{proposition}\label{path_rho1}
Let $\rho \co \pi_1(\O)\simeq \pi_1(\O_L)\ast_{\lan [\gamma] \ran}\pi_1(\O_R) \to SL(n,\Q)$ be a representation for which $\rho([\gamma])$ has $n$ distinct positive eigenvalues. 
Then there exists a path of representations $\rho_t \co \pi_1(\O) \to SL(n,\R)$ with $t\geq0$ such that
\begin{enumerate}
	\item $\rho_0 = \rho$,
	\item $\rho_t(\pi_1(\O)) = \lan \rho(\pi_1(\O_L)),\delta_t \rho(\pi_1(\O_R)) \delta_t^{-1}\ran$ for some $\delta_t\in SL(n,\R)$ which commutes with $\rho([\gamma])$, and 
	\item $\rho_t$ has image in $SL(n,\Q)$ for every $t\in \Q$.
\end{enumerate}
\end{proposition}
\proof The matrix $\rho([\gamma])$ is conjugate to a diagonal matrix $D$ with entries  $\lambda_1, \ldots, \lambda_n >0$ along its diagonal. 
Now for every $t>0$ define
\begin{equation}\label{deltat}
 	\delta_t = (t\rho([\gamma]) + I)\det(t\rho([\gamma]) + I)^{-\frac{1}{n}}
\end{equation}
Notice that $\det(t\rho([\gamma]) + I) = \det(tD + I) = \Pi_{k=1}^{n}(t\lambda_i +1)>0$,
so $t\rho([\gamma]) +I$ is invertible for all $t$. 
Then each $\delta_t$ is in $SL(n,\R)$ and we can check that $\delta_t$ commutes with $\rho([\gamma])$. 
Since $\rho$ is a rational representation, 
whenever $t\in \Q$ the matrix $t\rho([\gamma]) + I$ has rational entries and non-zero determinant. 

Let $\rho_t \co \pi_1(\O) \to SL(n,\R)$ be the representation such that $\rho_t(\pi_1(\O)) = \lan \rho(\pi_1(\O_L)),\delta_t \rho(\pi_1(\O_R)) \delta_t^{-1}\ran$.
Notice that $\rho_0 = \rho$ and that for every $t\in\Q$ the representation $\rho_t$ has image in $SL(n,\Q)$.

\QED

\subsection{Discarding Zariski closures}

For the rest of section \ref{sectionBending} we focus on the case where $n=2k+1$ is odd. 
Recall that in this case $SL(n,\R) \equiv PSL(n,\R)$. 

\begin{lemma}\label{Junique}
Let $\rho\co\Gamma \to SL(n,\R)$ be an irreducible representation and suppose there is a quadratic form $J$ such that $\rho(\Gamma )\subset SO(J)$. 
Then $J$ is unique up to scaling. 
\end{lemma}
\proof 
Suppose $\rho(\Gamma)< SO(J_1)\cap SO(J_2)$. Then for any $\rho(\gamma) \in \rho(\Gamma)$ we have that
	$$J_1^{-1} \rho(\gamma)J_1 = \rho(\gamma)^{-T} = J_2^{-1}\rho(\gamma)J_2,$$
which implies that 
	$\rho(\gamma)J_1J_2^{-1} = J_1J_2^{-1}\rho(\gamma).$
Since $n$ is odd, $J_1J_2^{-1}$ has a real eigenvalue $\lambda$.
Then $\Ker(J_1J_2^{-1} - \lambda I)$ is a non-zero invariant subspace for the irreducible representation $\rho$, which implies $J_1 = \lambda J_2$. 

\QED

\begin{proposition}\label{notZclosure}
Let  $\rho \co \pi_1(\O)\simeq \pi_1(\O_L)\ast_{\lan [\gamma] \ran}\pi_1(\O_R) \to SL(n,\R)$ be a representation in which  the restrictions $\rho|_{\pi_1(\O_L)}$ and $\rho|_{\pi_1(\O_R)}$ are  irreducible and $\rho([\gamma])$ has $n$ positive distinct eigenvalues. 
Suppose there is a quadratic form $J$ such that $\rho(\pi_1(\O))\subset SO(J)$. 
Then there exists a path of representations $\rho_t \co \pi_1(\O) \to SL(n,\R)$ such that
\begin{enumerate}
\item $\rho_0 = \rho$ and 
\item for each $t>0$ there is no quadratic form $\tilde{J}$ such that $\rho_t(\pi_1(\O))\subset SO(\tilde{J})$.
\end{enumerate}

\end{proposition}
\proof By proposition \ref{path_rho1} there are $\delta_t\in SL(n,\R)$ that commute with $\rho([\gamma])$, with which we can construct a path of representations $\rho_t\co \pi_1(\O)\to SL(n,\R)$ such that $\rho_0=\rho$ and $\rho_t(\pi_1(\O)) = \lan \rho(\pi_1(\O_L)),\delta_t \rho(\pi_1(\O_R)) \delta_t^{-1}\ran$.

Now fix $t>0$. Suppose there exists a quadratic form $\tilde{J}$ such that $\rho_t(\pi_1(\O)) \subset SO(\tilde{J})$. 
Since $\rho(\pi_1(\O)) \subset SO(J)$, in particular $\rho_t(\pi_1(\O_L)) = \rho_0(\pi_1(\O_L)) \subset SO(J)\cap SO(\tilde{J})$. 
The restriction $\rho_t|_{\pi_1(\O_L)}$ is  irreducible, so by lemma \ref{Junique} $J$ is a real multiple of $\tilde{J}$. 
Similarly, by construction $\rho_t (\pi_1(\O_R)) \subset SO(\delta_t J \delta_t^T)\cap SO(\tilde{J})$ and $\rho_t|_{\pi_1(\O_R)}$ is  irreducible too. 
Thus $\delta_t J \delta_t^T$ is also a multiple of $\tilde{J}$. 
This implies there is a $\lambda \in \R$ such that $\lambda J = \delta_t J \delta_t^T$ and then $\lambda^n = \det(\delta_t)^2 =1.$ 
Since $n$ is odd it must be that $\lambda =1 $ and we obtain $\delta_t \in SO(J)$. 
Given that
$$
(t\rho([\gamma]) + I) J (t\rho([\gamma])^T + I) 
					\ =\  t^2J \ +\  tJ(\rho([\gamma])^{T})^{-1} \ +\  tJ\rho([\gamma])^T\ +\ J, 
$$
having $J = \delta_tJ\delta_t^T$ would imply that  $ \mu I =  \rho([\gamma])^{-1} + \rho([\gamma])$ for some $\mu\in \R$.  
Recall that $\rho([\gamma])$ is conjugate to a diagonal matrix $D$ whose eigenvalues are all distinct. 
If  $\mu I = \rho([\gamma])^{-1} + \rho([\gamma])$ then by conjugating we would obtain that $\mu I = D^{-1} + D$, which is not the case given that $n>2$. 

\QED

\subsection{Representations of surface groups} 
Recall we are assuming that $\O$ is a 2-dimensional orientable connected closed orbifold of negative orbifold Euler characteristic. 
Such orbifolds are always finitely covered by a surface $S$ of genus greater than one, so  $\pi_1(S)$ is a finite index subgroup of $\pi_1(\O)$. 
Given a representation $\rho \co \pi_1(\O) \to G$ we will denote the restriction of $\rho$ to $\pi_1(S)$ by $\rho^S$. 

\begin{theorem}\label{Zdense_surfacegrps}
Suppose $\pi_1(\O)\simeq \pi_1(\O_L)\ast_{\lan [\gamma] \ran}\pi_1(\O_R)$ with $[\gamma]$ an infinite order element. 
Let $\rho \co  \pi_1(\O)  \to SL(n,\R)$ be an orbifold Fuchsian representation such that the restrictions $\rho|_{\pi_1(\O_L)}$ and $\rho|_{\pi_1(\O_R)}$ are irreducible.
If $S$ is a surface finitely covering $\O$ then there exists a path of representations $\rho^S_t \co \pi_1(S) \to SL(n,\R)$ such that $\rho^S_0 = \rho^S$ and $\rho_t^S$ is a Zariski dense surface Hitchin representation for each $t>0$. 
\end{theorem}
\proof 
Since $\rho\co \pi_1(\O) \to SL(n,\R)$ is an orbifold Hitchin representation with odd $n=2k+1$ and $[\gamma]$ has infinite order, then $\rho([\gamma])$ has $n$ positive distinct real eigenvalues. 
Moreover, since $\rho$ is Fuchsian its image is contained in a conjugate of $SO(k,k+1)$. 
Using proposition \ref{notZclosure} we obtain a path of representations $\rho_t \co \pi_1(\O) \to SL(n,\R)$ such that  $\rho_0 = \rho$ and for each $t>0$ there is no real quadratic form $J$ such that $\rho_t(\pi_1(\O))\subset SO(J)$. 
By proposition \ref{Zdense} each $\rho_t(\pi_1(S))$ is Zariski dense in $SL(n,\R)$.

Now consider the continuous path $[\rho_t] \in \text{Rep}(\pi_1(\O),PGL(n,\R))$ for $t\geq0$. 
Its image is connected so all $PGL(n,\R)$-conjugacy classes $[\rho_t]$ are contained in the same connected component of $\text{Rep}(\pi_1(\O),PGL(n,\R))$. 
Because the representation $\rho_0 = \rho$ is Fuchsian, $[\rho_0]$ is in the Hitchin component $\text{Hit}(\pi_1(\O), PGL(n,\R))$ and so is every $[\rho_t]$. 
Thus, by theorem \ref{orbifoldHitchin}, each $\rho_t$ is discrete, faithful and strongly irreducible. 
Since $\pi_1(S)$ has finite index in $\pi_1(\O)$, each restriction $\rho^S_t \co \pi_1(S) \to SL(n,\R)$ is irreducible. 
In particular $\rho_0^S$ is a surface Fuchsian representation. 
Then $[\rho^S_t]$ is a continuous path in $\text{Rep}^+(\pi_1(S),SL(n,\R))$ with $[\rho_0^S] \in \text{Hit}(\pi_1(S),SL(n,\R))$. 
Since the Hitchin component is path connected $[\rho^S_t] \in \text{Hit}(\pi_1(S),SL(n,\R))$ for all $t\geq0$. 

\QED

To finish this section notice that the construction of the path of Zariski dense representations in the previous theorem is based on proposition \ref{path_rho1}, so we may add the assumption of $\rho(\pi_1(\O))\subset SL(n,\Q)$ to obtain that the image of every $\rho_t$ is in $SL(n,\Q)$ for every $t\in \Q$.

\begin{corollary}\label{rational_reps}
Let $\rho \co \pi_1(\O)\to PSL(n,\Q)$ be a representation satisfying the assumptions of theorem \ref{Zdense_surfacegrps}. 
If $S$ is a surface finitely covering $\O$ then there exists a path $\rho^S_t \co \pi_1(S) \to SL(n,\R)$ of Hitchin representations such that $\rho^S_0 = \rho^S$, $\rho_t^S$ is Zariski dense  for each $t>0$ and $\rho_t^S$ has image in $SL(n,\Q)$ for every $t\in \Q$. 
\end{corollary}

\section{Representations of $\pi_1(\O_{3,3,3,3})$}\label{sectionO3333}

In this section we look at the orbifold $\O_{3,3,3,3}$ and find a Fuchsian representation $\rho \co \pi_1(\O_{3,3,3,3})\to SL(n,\Z)$ satisfying the assumptions of corollary \ref{rational_reps}. 

\subsection{The orbifold $\O_{3,3,3,3}$} 

In what follows we focus on the triangle group $\Delta(3,4,4) \subset PSL(2,\R)$. 
If we let $T$ be the hyperbolic triangle with angles $\{\frac{\pi}{3}$, $\frac{\pi}{4}$, $\frac{\pi}{4}\}$, then the 
 generators of $\Delta(3,4,4)$ are the rotations $x$ and $y$ by $\frac{2\pi}{3}$ and $\frac{\pi}{2}$ around the corresponding vertices of $T$.
This group has presentation
\begin{equation}\label{triangle_presentation}
	\Delta(3,4,4) = \lan x, y  \ | \ x^3 = y^4 = (xy)^4 =  1 \ran.
\end{equation}
The fundamental domain for the action of $\Delta(3,4,4)$ on $\H$ is a quadrilateral with angles
 $\{\frac{\pi}{2},\frac{\pi}{3},\frac{\pi}{2},\frac{\pi}{3}\}$. 
 The quotient $\H/\Delta(3,4,4)$ is homeomorphic to the orbifold $S^2(3,4,4)$ whose underlying topological space is $S^2$ and has three cone points of orders $3,\ 4$ and $4$. 
This defines, up to conjugation, an isomorphism $\pi_1(S^2(3,4,4))\to \Delta(3,4,4)\subset PSL(2,\R)$.

\begin{figure}[H]
\begin{center}
\includegraphics[scale=0.3]{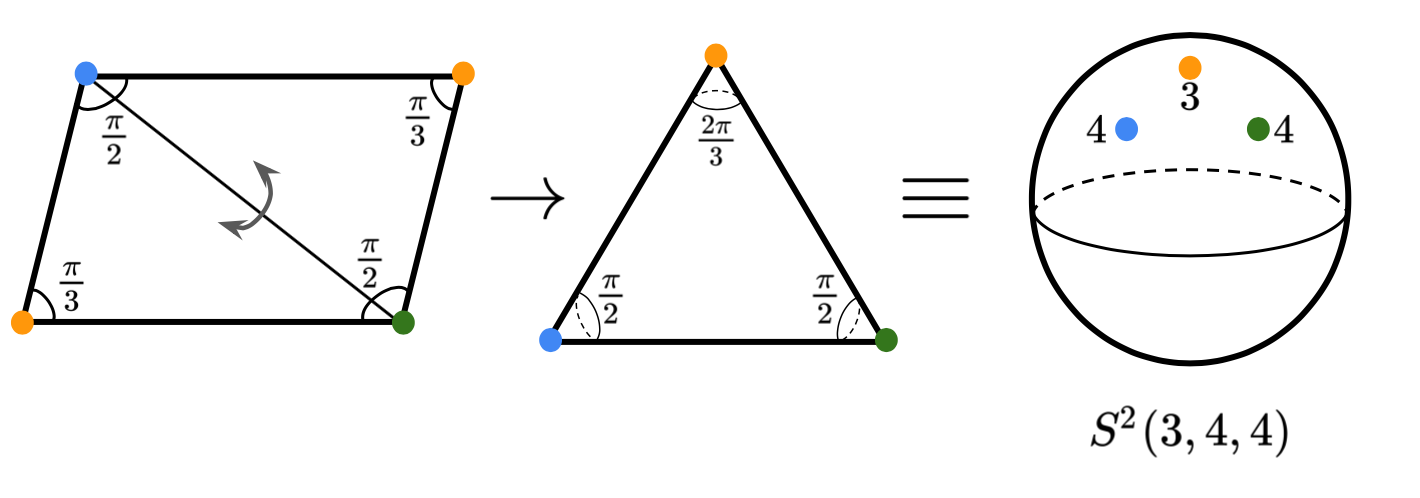} \ \ \ \ \ \ 
\caption{Orbifold $S^2(3,4,4)$}
\end{center}
\label{fig:flagsapt3}
\end{figure}

Let
	$\theta_1 = x $ and $ \theta_i = y\theta_{i-1}y^{-1} $ for $ i = 2,3,4$, 
then $\lan \theta_1, \ldots, \theta_4 \ran$ the quotient of $\H$ by the action of $\lan \theta_1, \ldots, \theta_4 \ran$ is homeomorphic to the orbifold $\mathcal{O}_{3,3,3,3}$ with underlying topological space $S^2$ and  4 cone points of order 3. 
By construction, we obtain that $\mathcal{O}_{3,3,3,3}$ is an index four orbifold covering of $S^2(3,4,4)$. 
If $ \gamma_1, \ldots, \gamma_4$ are loops around the cone points of $\mathcal{O}_{3,3,3,3}$, then the orbifold fundamental group has the presentation
$$
	\pi_1(\mathcal{O}_{3,3,3,3}) = \lan \gamma_1, \ldots, \gamma_4 \ | \  \gamma_1^3 = \ldots = \gamma_4^3 =   \gamma_1  \gamma_2  \gamma_3  \gamma_4 = 1\ran.
$$
Identifying each $ \gamma_i$ with the rotation $\theta_i$ gives an isomorphism $\pi_1(\mathcal{O}_{3,3,3,3})\cong \lan \theta_1, \ldots, \theta_4\ran$ which defines (up to conjugation) a discrete and faithful representation 
\begin{equation}\label{init_rep}
	\sigma \co \pi_1(\mathcal{O}_{3,3,3,3}) \to \Delta(3,4,4) < PSL(2,\R).
\end{equation}

\begin{lemma}\label{sigma_Zdense}
The representation $\sigma\co\pi_1(\O_{3,3,3,3})\to PSL(2,\R)$ defined in (\ref{init_rep}) is Zariski dense.
\end{lemma}
\proof 
We will check that the group $ \sigma(\pi_1(\O_{3,3,3,})) = \lan \theta_1, \ldots, \theta_4 \ran < \Delta(3,4,4)$ is Zariski dense. 
Hyperbolic triangles with the same angles are isometric, 
so we can fix the hyperbolic triangle with angles $\{\frac{\pi}{3}, \frac{\pi}{4},\frac{\pi}{4}\}$ by placing it symmetrically along the $y$-axis in the upper-half plane. 
By having the generators $x,y$ of $\Delta(3,4,4)$ defined in (\ref{triangle_presentation}) in rational canonical form we obtain that:
\begin{eqnarray}\label{gens}
	x = \begin{bmatrix} 0 & -1 \\ 1 & 1 \end{bmatrix} \  \text{ and } \ 
	y = \begin{bmatrix} 0 & -1 - \sqrt{2} \\ -1 +\sqrt{2} & \sqrt{2} \end{bmatrix}.
\end{eqnarray}
This choice of generators fixes a representative in the conjugacy class of the representation $\sigma$. 
Notice that $\theta_2\theta_1 = yxy^{-1}x$ is an infinite order element in $\Delta(3,4,4)$ and is therefore hyperbolic.
By using the matrices in (\ref{gens}) we can explicitly find $P,D \in PGL(n,\R)$ with $D$ diagonal so that $P^{-1}(\theta_2\theta_1)P = D$. 
It suffices then to see that the conjugated representation $P^{-1}\sigma P$ is Zariski dense. 
Let $H$ be the Zariski closure of  $P^{-1}\sigma(\pi_1(\mathcal{O}_{3,3,3,3}))P$ in $PSL(2,\R)$ and $\mathfrak{h}$ its Lie algebra.  
First notice that the Zariski closure of $\lan D \ran$ is the algebraic torus 
whose Lie algebra is the span of $X_1 = \begin{pmatrix} 1 & 0 \\ 0 & -1 \end{pmatrix}$. 
Taking 
$ X_2 = \text{Ad}_{P^{-1}\theta_1\theta_2P}(X_1)$ and $X_3 = \text{Ad}_{P^{-1}\theta_1^2\theta_2P}(X_1)$ we obtain three linearly independent vectors in  $\mathfrak{h}$. Then $\dim(\mathfrak{h})=3 = \dim(\mathfrak{sl}(2,\R))$ so the two algebras must coincide and so $H=PSL(2,\R)$. 

\QED

\subsection{Rational representations of $\pi_1(\O_{3,3,3,3})$}

We will now focus on the case $n=2k+1$ and the representation
$
	\omega_n \circ \sigma \co  \pi_1(\mathcal{O}_{3,3,3,3}) \to SL(n,\R),
$
where $\sigma$ is the representation defined in (\ref{init_rep}) and $\omega_n\co PSL(2,\R) \to PSL(n,\R) =SL(n,\R)$ the irreducible representation introduced in \ref{subsec_irrep}. 
Since $\omega_n\circ \sigma$ is an orbifold Fuchsian representation, it is irreducible. 
The following result implies that we can conjugate $\omega_n\circ \sigma$ to obtain an integral representation 
\begin{eqnarray}\label{rho_0}
	\rho \co \pi_1(\mathcal{O}_{3,3,3,3}) \to SL(n,\Z) < SL(n,\R). 
\end{eqnarray}

\begin{proposition}[\cite{Long_This_20_SLodd} thm. 2.1 ]\label{conjugateZ}
Let $\omega_n\co PSL(2,\R) \to PSL(n,\R)$ be the unique irreducible representation between these groups. 
Then for every odd $n$ the restriction $\phi_n = \omega_n |_{\Delta(3,4,4)}$ is conjugate to a representation $\rho_n \co \Delta(3,4,4) \to PSL(n,\Z)$. 
\end{proposition}

Now let $\gamma \subset \O_{3,3,3,3}$ be a simple closed loop dividing $\O_{3,3,3,3}$ into two orbifolds $\O_L$ and $\O_R$ which share $\gamma$ as their common boundary and have two cone points of order 3 each. 
Then $[\gamma] \in \pi_1(\O_{3,3,3,3})$ is an infinite order element and 
	$\pi_1(\O_{3,3,3,3}) \simeq \pi_1(\O_L)\ast_{\lan [\gamma] \ran}\pi_1(\O_R).$

\begin{proposition}\label{t_abs_irred}
Let $\rho\co \pi_1(\O_{3,3,3,3})  \simeq \pi_1(\O_L)\ast_{\lan [\gamma] \ran}\pi_1(\O_R) \to PSL(n,\Z)$ be the representation defined in (\ref{rho_0}). Then the restrictions of $\rho$ to $\pi_1(\O_L)$ and $\pi_1(\O_R)$ are  irreducible.
\end{proposition}
\proof
To see that $\rho|_{\pi_1(\O_L)}$ is  irreducible it suffices to see that the restriction of $\omega_n\circ \sigma$ to $\pi_1(\O_L)$ is  irreducible. 
By the proof of lemma \ref{sigma_Zdense} we have that $\sigma(\pi_1(\O_L))$ is Zariski dense in $PSL(2,\R)$. 
To see that the representation $\omega_n\co\sigma(\pi_1(\O_L)) \to PSL(n,\R)$ is irreducible, it is enough to check that the Zariski closure of its image is irreducible. 
This holds since $\omega_n \co PSL(2,\R) \to PSL(n,\R)$ is an  irreducible representation and a morphism of algebraic groups, so
$\omega_n(PSL(2,\R)) = \omega_n(\overline{\sigma(\pi_1(\O_L)}) \subseteq \overline{\omega_n\circ\sigma(\pi_1(\O_L))}$.

To see $\rho|_{\pi_1(\O_R)}$ is irreducible it is enough to notice that the proof of \ref{sigma_Zdense} also holds for $\pi_1(\O_R)$ by using the generators $\theta_3$ and $\theta_3$ instead of $\theta_1$ and $\theta_2$. 

\QED

Knowing that $\rho$ is an integral orbifold Fuchsian representation, the previous proposition shows $\rho$ satisfies the assumptions of theorem \ref{Zdense_surfacegrps}. 
Thus we obtain the following application of corollary \ref{rational_reps}.

\begin{theorem}
For every surface $S$ finitely covering the orbifold $O_{3,3,3,3}$ and every odd $n>1$ there exists a path of Hitchin representations $\rho_t \co \pi_1(S) \to SL(n,\R)$, so that
\begin{enumerate}
	\item  $\rho_0(\pi_1(S)) \subset SL(n,\Z)$, 
	\item $\rho_t$ is Zariski dense for every $t>0$ and
	\item $\rho_t(\pi_1(S)) \subset SL(n,\Q)$ for every $t\in \Q$. 
\end{enumerate}
\end{theorem}

\bibliographystyle{ieeetr}
\bibliography{SUBMISSION_z_dense}

\end{document}